\documentclass[12pt, reqno]{amsart}
\usepackage{amsmath,amsfonts,amsthm,amsopn,amssymb}
\usepackage{cite,marginnote}
\usepackage{color,enumitem,graphicx}
\usepackage[colorlinks=true,urlcolor=blue,
citecolor=red,linkcolor=blue,linktocpage,pdfpagelabels,
bookmarksnumbered,bookmarksopen]{hyperref}
\usepackage[english]{babel}

\usepackage[left=2.9cm,right=2.9cm,top=2.8cm,bottom=2.8cm]{geometry}
\usepackage[hyperpageref]{backref}

\numberwithin{equation}{section}

\makeindex
\newtheorem{theorem}{Theorem}[section]
\newtheorem{lemma}{Lemma}[section]
\newtheorem{corollary}{Corollary}[section]
\newtheorem{proposition}{Proposition}[section]
\newtheorem{remark}{Remark}[section]
\newtheorem{example}{Example}[section]
\newtheorem{definition}{Definition}[section]
\newtheorem{theoremletter}{Theorem}

\newcommand{\ud}{\mathrm{d}}
\newcommand{\RN}{\mathbb R^N}
\newcommand{\om}{\Omega}

\newcommand{\iy}{\infty}

\newcommand{\s}{\section}

\newcommand{\DD}{\Delta}
\newcommand{\g}{\gamma}
\newcommand{\G}{\Gamma}
\newcommand{\na}{\nabla}

\newcommand{\la}{\lambda}

\newcommand{\pa}{\partial}
\newcommand{\si}{\sigma}

\newcommand{\R}{\mathbb R}
\newcommand{\al}{\alpha}

\newcommand{\ti}{\tilde}
\newcommand{\wti}{\widetilde}

\newcommand{\re}[1]{(\ref{#1})}

\newcommand{\rg}{\rightarrow}

\newcommand{\e}{\varepsilon}

\newcommand{\lab}{\label}
\newcommand{\bl}{\begin{lemma}}
\newcommand{\el}{\end{lemma}}
\newcommand{\bd}{\begin{definition}}
\newcommand{\ed}{\end{definition}}
\newcommand{\bc}{\begin{corollary}}
\newcommand{\ec}{\end{corollary}}
\newcommand{\bp}{\begin{proof}}
\newcommand{\ep}{\end{proof}}
\newcommand{\bx}{\begin{example}}
\newcommand{\ex}{\end{example}}
\newcommand{\bi}{\begin{exercise}}
\newcommand{\ei}{\end{exercise}}
\newcommand{\bo}{\begin{proposition}}
\newcommand{\eo}{\end{proposition}}
\newcommand{\br}{\begin{remark}}
\newcommand{\er}{\end{remark}}
\newcommand{\be}{\begin{equation}}
\newcommand{\ee}{\end{equation}}
\newcommand{\ba}{\begin{align}}
\newcommand{\ea}{\end{align}}
\newcommand{\bn}{\begin{enumerate}}
\newcommand{\en}{\end{enumerate}}
\newcommand{\bg}{\begin{align*}}
\newcommand{\bcs}{\begin{cases}}
\newcommand{\ecs}{\end{cases}}

\newcommand{\bean}{\begin{eqnarray*}}
\newcommand{\eean}{\end{eqnarray*}}

\title[Existence and concentration of positive solutions]{Existence and concentration of positive solutions for nonlinear Kirchhoff type problems with a general critical nonlinearity}

\author[J. Zhang]{Jianjun Zhang}
\author[D. G.\ Costa]{David G.\ Costa}
\author[J. M.\ do \'O]{Jo\~ao Marcos do \'O}

\address[J. J.\ Zhang]{\newline\indent College of Mathematics and Statistics
\newline\indent
Chongqing Jiaotong University
\newline\indent
Chongqing 400074, PR China
\newline\indent and
\newline\indent Dip. di Scienza e Alta Tecnologia
\newline\indent
Universit\`{a} degli Studi dell'Insubria
\newline\indent
via Valleggio 11, 22100 Como,Italy}
\email{\href{mailto:zhangjianjun09@tsinghua.org.cn}{zhangjianjun09@tsinghua.org.cn}}

\address[D. G.\ Costa]{\newline\indent Department of Mathematics Sciences
\newline \indent
University of Nevada Las Vegas
\newline\indent Las Vegas, P.O. Box No. 454020, NV, USA}
\email{\href{mailto:costa@unlv.nevada.edu}{costa@unlv.nevada.edu}}

\address[J. M. do \'O]{\newline\indent Department of Mathematics
\newline\indent
Federal University of Para\'{\i}ba
\newline\indent
58051-900, Jo\~ao Pessoa-PB, Brazil}
\email{\href{mailto:jmbo@pq.cnpq.br}{jmbo@pq.cnpq.br}}

\thanks{Corresponding author: J.~M.~do~\'{O}. \\
Research partially supported by The National Institute of Science and Technology of Mathematics ICNT-Mat,
CAPES and CNPq/Brazil, J.\ Zhang was partially supported by the Science Foundation of Chongqing Jiaotong University(15JDKJC-B033)}


\begin{document}

\begin{abstract}
We are concerned with the following Kirchhoff type equation
$$
-\e^2M\left(\e^{2-N}\int_{\RN}|\na u|^2\, \ud x\right)\Delta u+V(x)u=f(u),\ x\in\RN,\ \ N\ge2,
$$
where $M\in C(\R^+,\R^+)$, $V\in C(\RN,\R^+)$ and $f(s)$ is of critical growth. In this paper, we construct a localized bound state solution concentrating at a local minimum of $V$ as $\e\rg 0$ under certain conditions on $f(s)$, $M$ and $V$. In particular, the monotonicity of $f(s)/s$ and the Ambrosetti-Rabinowitz condition are not required. 	
\end{abstract}
\maketitle

\section{Introduction and main results}

\noindent

In this paper, we are concerned with existence and concentration of positive solutions to the following Kirchhoff type equations
\be\lab{q2}
-\e^2M\left(\e^{2-N}\int_{\RN}|\na v|^2\, \ud x\right)\Delta v+V(x)v=f(v)\ \mbox{in}\ \RN,
\ee
where $M\in C(\R^+,\R^+)$, $V\in C(\RN,\R^+)$, $N\ge2$ and $\e>0$. In the sequel, we assume that
the potential $V$ satisfies
\begin{itemize}

\item[$(V1)$] $V\in C(\RN,\R)$ and $0<V_0=\inf_{x\in\RN}V(x)$,

\item[$(V2)$] there is a bounded domain $O\subset \RN$ such that
$$
m\equiv\inf_{x\in O}V(x)<\min_{x\in \partial O}V(x),
$$
\end{itemize}
and $M$ satisfies $(M1)$ if $N=2$ and $(M1)$-$(M5)$ if $N\ge3$ below:
\begin{itemize}

\item[$(M1)$] $0<m_0=\inf_{t\in\R^+}M(t)$;

\item[$(M2)$] if $\hat{M}(t)=\int_0^tM(s)\,\ud s$, we have
$$
\liminf_{t\rg+\iy}\left[\hat{M}(t)-(1-\frac{2}{N})tM(t)\right]=+\iy;
$$

\item[$(M3)$] $\lim_{t\rg+\iy}M(t)/t^{2/(N-2)}=0$;

\item[$(M4)$] $M(t)$ is nondecreasing for $t\in\R^+$;

\item[$(M5)$] $M(t)/t^{2/(N-2)}$ is nonincreasing for $t\in\R^+$.
\end{itemize}
In 2014, G. M. Figueiredo et al. \cite{Figu} considered the concentration phenomenon of the above problem \re{q2} in the subcritical case. The authors assumed that
\begin{itemize}

\item[$({f_0})$] $f\in C(\R, \R)$ and $f(t)=0$ if $t\le0$;

\item[$({f_1})$]  $-\iy<\liminf_{t\rg0}f(t)/t\le\limsup_{t\to 0}f(t)/t<V_0;$

\item[$({f_2})$]  $\lim_{t\rg\iy}f(t)/e^{\al t^2}=0$ for any $\al>0$ if $N=2$, $\lim_{t\rg\iy}f(t)/t^{(N+2)/(N-2)}=0$ if $N\ge3$;

\item[$({f_3})$]  there exists $\xi>0$ such that $\xi^2{m}<2F(\xi)$, where  $F(s):=\int_0^sf(t) \, \ud t$.
\end{itemize}

\noindent  Let
$$
\mathcal{M}\equiv\{x\in O: V(x)=m\}.
$$

\begin{theoremletter}\label{thmA}
{\rm(see \cite{Figu})} {\it Assume $(V1)$-$(V2)$, $(f_0)$-$(f_3)$, $(M1)$ if $N=2$ and $(M1)$-$(M5)$ if $N\ge3$. Then, for sufficiently small $\e>0$, \eqref{q2} admits a positive solution $v_{\e}$ which satisfies:
\begin{itemize}
\item[(i)] there exists a maximum point $x_\e$ of $v_\e$ such that $\lim_{\e\rg 0}dist(x_\e,\mathcal{M})=0$ and, for any such $x_\e$, $w_\e(x)\equiv v_\e(\e x+x_\e)$ converges (up to a subsequence) uniformly to a least energy solution of
\begin{eqnarray*} -M(\|\na u\|_2^2)\DD u+mu=f(u),\ \ u>0,\ \ u\in H^1(\RN); \end{eqnarray*}
\item[(ii)] $v_\e(x)\le C\exp(-\frac{c}{\e}|x-x_\e|)$ for some $c,C>0$.
\end{itemize}}
\end{theoremletter}

Before stating our main result, we shall introduce the main hypotheses on $f$. In what follows, we assume that $f\in
C(\R^+,\R^+)$ and satisfies
\begin{itemize}
\item[(F1)] $\lim_{t\rg 0^+}{f(t)}/{t}=0$.
\item[(F2)] If $N=2$, $\lim_{t\rg+\iy}{f(t)}/{\exp(\al t^2)}=\left\{
\begin{array}{ll}
0,\ \ \ \ \ \ \forall \al>4\pi,\\
+\iy,\ \ \ \forall \al<4\pi.
\end{array}
\right.
$

\noindent If $N\ge3$, $\lim_{t\rg\iy}f(t)/t^{(N+2)/(N-2)}=1$.
\vskip0.1in
\item [(F3)] If $N=2$, there exists $\beta_0>e\, m/(2\pi)$ such that $$\lim\limits_{|t|\rg+\iy}tf(t)/\exp{(4\pi t^2)}\ge\beta_0.$$
 \vskip0.1in
\noindent If $N\ge 3$, there exist $\la>0$ and $p<2^\ast$ such that $$f(t)\ge t^{(N+2)/(N-2)}+\la t^{p-1}, \ t\ge0,$$ where $p$ and $\la$ satisfy {\bf one} of the following conditions:
\begin{itemize}
\item[(i)] $p\in(2,2^\ast)$ and $\la>0$ if $N\ge4$;
\item[(ii)] $p\in(4,2^\ast)$ and $\la>0$ if $N=3$;
\item[(iii)] $p\in(2,4]$ and $\la>0$ large enough if $N=3$.
\end{itemize}
\end{itemize}

The main theorem of this paper reads as

\begin{theorem}
\lab{Theorem 1} {\it Assume $(V1)$-$(V2)$, $(F1)$-$(F3)$, $(M1)$ if $N=2$ and $(M1)$-$(M5)$ if $N\ge3$. Then, for sufficiently small $\e>0$, \eqref{q2} admits a positive solution $v_{\e}$, which satisfies
\begin{itemize}
\item[(i)] there exists a maximum point $x_\e$ of $v_\e$ such that $\lim_{\e\rg 0}dist(x_\e,\mathcal{M})=0$ and for any such $x_\e$, $w_\e(x)\equiv v_\e(\e x+x_\e)$ converges (up to a subsequence) uniformly to a least energy solution of
\begin{eqnarray*} -M(\|\na u\|_2^2)\DD u+mu=f(u),\ \ u>0,\ \ u\in H^1(\RN); \end{eqnarray*}
\item[(ii)] $v_\e(x)\le C\exp(-\frac{c}{\e}|x-x_\e|)$ for some $c,C>0$.
\end{itemize}}
\end{theorem}

Now, let us give some more background for \eqref{q2}. For $\e=1$ and a bounded domain $\Omega\subset\RN$, \re{q2} is reduced to
\be\lab{qq2}
-M\left(\int_{\Omega}|\na v|^2\right)\DD v+V(x)v=f(v),\,\,x\in\Omega.
\ee
Equation \re{qq2} arises when one seeks steady states to the time-dependent wave type equation
\be\lab{K}
u_{tt}-M\left(\int_{\Omega}|\na v|^2\right)\DD v+V(x)v=f(v),\,\,(t,x)\in(0,T)\times\Omega,
\ee
as well as when looking for the standing wave to the time-dependent Schr\"odinger equation
$$
i\hbar\frac{\pa\Psi}{\pa t}=-\frac{\hbar^2}{2m}\DD \Psi+V(x)\Psi-f(\Psi),\ \ (t,x)\in\R\times\R^N,
$$
Problem \re{K} was proposed by Kirchhoff in \cite{K} with $M(t)=a+bt$ and $N=1$. After the works of  Kirchhoff \cite{K} and Lions \cite{Lions}, the Kirchhoff problem \re{qq2} have been paid much attention. For more background, we refer to \cite{Figu} and the references therein.
\vskip0.1in
For the case $M(t)=1$, Problem \re{q2} reads
\be\lab{qq1}
-\e^2\DD v+V(x)v=f(v),\ \ v\in H^1(\RN).
\ee
In the last decades, considerable attention has been paid to problem \eqref{qq1}. An interesting class of solutions of
\eqref{qq1} consists of families of solutions which develop a spike shape around some point in $\RN$ as $\e\rg 0$. From the physical point
of view, these solutions are referred to as semiclassical states, as they describe the transition from classical mechanics to quantum mechanics.

\vskip0.1in
After the celebrated work of Floer and Weinstein \cite{F-W}, Problem \re{qq1} has been studied by many researchers. Here we only refer to \cite{Oh,Rab,WX,Felmer} and the references therein. But in these works, the nonlinearity $f(s)$ is basically required to satisfy the $Ambrosetti$-$Rabinowitz$ condition:
\be\tag{$A$-$R$}
0<\mu\int_{R}f(s)\,\ud s\le sf(s),\mu>2\,\,\,\mbox{for}\,\,s\not=0,
\ee
and a monotonicity condition:
\be\tag{$M$}
f(s)/|s|\,\,\,\mbox{is strictly increasing for}\,\,s\not=0.
\ee
A natural question is whether these results hold for more general nonlinear terms $f(s)$,
particularly, without $(A$-$R)$ and the monotonicity condition $(M)$. In 2007, J. Byeon and L. Jeanjean \cite{byeon} gave a positive answer for $N\ge3$. Precisely, they assume $(V1)$-$(V2)$ as in Problem \re{q2}. Then under the Berestycki-Lions conditions $(F1)$-$(F3)$ in \cite{byeon}, they constructed a spike solution for \eqref{qq1} around the local minimum of $V$ stated in $(V2)$. In 2008, Byeon, Jeanjean and Tanaka \cite{BJT} used a similar argument to \cite{byeon} to obtain a corresponding result for \eqref{q2} in the cases $N=1,2$. Moreover, the hypotheses in \cite{byeon,BJT} are almost optimal.
For the critical case with general nonlinearities, we refer to the recent works \cite{zhang-chen-zou,JZ1}. Through all these works above, the assumption $\inf_{x\in\RN}V(x)>0$ was imposed. It is easy to see that if $\inf_{x\in\RN}V(x)<0$, there
exists no solution of \re{qq1} for small $\e>0$. Thus, the case $\inf_{x\in\RN}V(x)=0$ is called the critical frequency. In \cite{byeon-wang} Byeon and Wang
gave the breakthrough for this condition. If $\liminf_{|x|\rg\iy}V(x)>0=\inf_{x\in\RN}V(x)$, Byeon and Wang \cite{byeon-wang} proved the existence of solutions concentrating on an isolated component of $\mathcal{Z}=\{x\in\RN:V(x)=0\}$. For further related result, we here also refer to Ambrosetti-Wang \cite{A-W}, Cao-Noussair \cite{Cao1}, Cao-Peng \cite{Cao2} and Cao-Noussair-Yan \cite{Cao3}, Moroz-Schaftingen \cite{Moroz} and the references therein.

\vskip0.1in

For the case $M(t)=a+bt$ and $N=3$, Problem \re{q2} reads
\be\lab{qq11}
-\left(\e^2 a+\e b\int_{\R^3}|\na v|^2\,\ud x\right)\DD v+V(x)v=f(v),\ \ v\in H^1(\R^3).
\ee
By the Nehari manifold method, X. He and W. Zou \cite{Zou} considered the existence and concentration of ground sate solutions to \re{qq11} in the subcritical case. Later, Wang et al.\cite{Tian} obtained similar results as in \cite{Zou} in the critical case. However, $(A$-$R)$ (with $\mu>3$) or the monotonicity condition
\be\tag{$M'$}
f(s)/s^3\,\,\,\mbox{is strictly increasing for}\,\,s>0
\ee
is required. Moreover, in \cite{Zou, Tian}, $f$ satisfies
$$
\lim_{s\rg0}f(s)/s^3=0.
$$
More recently, Y. He and G. Li \cite{HeLi} considered the existence and concentration of positive solutions to \re{qq11} with $f(s)=\la |s|^{p-2}s+s^5$. In \cite{HeLi}, with $p\in(2,4]$, the nonlinearity $f$ does not satisfies $(A$-$R)$ and the monotonicity condition $(M')$. Later, under the same assumptions on $f$ introduced in \cite{Zhang-Zou}, Y. He \cite{HH} extended the result in \cite{zhang-chen-zou} to the Kirchhoff problems. Here we should point out that in \cite{HeLi,HH,zhang-chen-zou}, the authors only considered the higher dimensional case ($N\ge3$) and the main ingredient used is indeed a Brezis-Nirenberg type argument. However, it seems very difficult to be adopted to deal with problem \re{q2} involving critical growth with respect to the Trudinger-Moser inequality. To the best of our knowledge, there are few results on the existence and concentration of solutions to \re{q2} involving a general critical nonlinearity in any dimension $N\ge2$. For the subcritical case, Figueiredo et al. in \cite{Figu} used similar arguments as in \cite{byeon,BJT} to get corresponding results for \re{q2}. Precisely, with the Berestycki-Lions conditions $(F1)$-$(F3)$ in \cite{byeon} or $(f_1)$-$(f_3)$ in \cite{BJT}, they obtained spike solutions around a local minimum of $V$.

It is natural to ask whether the result \cite{Figu} holds for more general nonlinear terms $f(s)$ in the critical case and for any dimension $N\ge2$. The main goal of this paper is two-fold. On one hand, we provide a new approach to deal with the \underline{critical case} for Kirchhoff-type problems in any dimension. The \underline{subcritical case} was considered by Figueiredo et al in \cite{Figu} as already pointed out (cf.~Theorem \ref{thmA}). Our approach when applied in the ``critical dimension'' $N=3$ is also considerably simpler than the one by Y. He and G. Li in \cite{HeLi} and Y. He \cite{HH}, in which the authors considered the Kirchhoff case $M(t)=a + bt$.

On the other hand, we also provide the concentration behavior of the corresponding ``semi-classical states'' $u_\e$, as $\e\to 0$. We point out that we allow critical perturbations which can locally go above the critical Sobolev exponent $2^{\star}:=\frac{2N}{N-2}$ (for $N\geq 3$) in the sense of assumption $(F3)$ in Theorem~\ref{thmA}, but in the corresponding critical situation. Our approach was inspired by the papers \cite{Figu,Azzollini} of Figueiredo et al and of Azzollini respectively, who provided a homeomorphism between the ground states of Kirchhoff equation and a related semilinear local elliptic equation.

The paper is organized as follows. Section 2 is devoted to the study of the so-called limit problem \re{q5} of \eqref{q2} (see below). The compactness of the set of ground sate solutions is proved. Section 3 is devoted to the proof of Theorem \ref{Theorem 1} by using the truncation approach in \cite{JZ1}.\vskip0.1in

\s{The limit problem}

\renewcommand{\theequation}{2.\arabic{equation}}

Since we are interested in the positive solutions of \re{q2}, from now on we may assume that $f(t)= 0$ for $t \le
0.$ In this case any weak solution of \re{q2} is positive by the maximum principle. The following equation when $m$ as in $(V2)$ is called the limiting equation of \eqref{q2}
\be\lab{q5}
-M(\|\na u\|_2^2)\DD u+mu=f(u),\ \ u\in H^1(\RN).
\ee
Define
$$
L_m(u)=\frac{1}{2}\hat{M}(\|\na u\|_2^2)+\frac{m}{2}\int_{\RN}u^2 \, \ud x-\int_{\RN}F(u)\, \, \ud x,\ \ u\in H^1(\RN)\,,
$$
and set
$$
\mathcal{N}_m:=\left\{u\in H^1(\RN)\setminus\{0\}: L_m'(u)=0, u(0)=\max_{x\in \RN}u(x)\right\}.
$$
Then we introduce the set $S_m$ of ground state solutions and the least energy $E_m$ of \eqref{q5} as follows:
$$
S_m:=\left\{u\in\mathcal{N}_m: L_m(u)=E_m\right\},\ \ \ E_m:=\inf_{u\in \mathcal{N}_m}L_m(u).
$$
Now, we give a result about $S_m$, whose proof follows from Lemma \ref{lemma1}-\ref{lemma3} below.
\bo\lab{prop1}\ Under the assumptions in Theorem \ref{Theorem 1} one has $S_m\not=\phi$. Moreover,
\begin{itemize}
\item [$(i)$] any $U\in S_m$ is such $U\in C^2(\RN)\cap L^\iy(\RN)$ and is radially symmetric;
\item [$(ii)$] $S_m$ is compact in $H^1(\RN)$;
\item [$(iii)$] $0<\inf\{\|U\|_{\infty}:U\in S_m\}\le \sup\{\|U\|_{\infty}:U\in S_m\}<\iy$;
\item [$(iv)$] there exist constants $C,c>0$ independent of $U\in S_m$ such that
$$
|D^{\al}U(x)|\le C\exp(-c|x|), \,x\in \RN\ \ \mbox{for\ \ $|\al|=0,1$.}
$$
\end{itemize}
\eo

We note that equation \re{q5} is nonlocal due to the presence of the term $M(\|\na u\|_2^2)$. Namely, \re{q5} is no longer a pointwise identity, which causes some mathematical difficulties in studying the properties of $S_m$. To overcome this difficulty, we use an idea introduced in \cite{Azzollini} and developed in \cite{Figu} to reduce equation.~\re{q5} to a local problem. Precisely, we consider
\be\lab{q51}
-\DD u+mu=f(u),\ \ u\in H^1(\RN),
\ee
whose energy functional is given by
$$
\ti{L}_m(u)=\frac{1}{2}\int_{\RN}\left(|\na u|^2+mu^2\right)\, \ud x-\int_{\RN}F(u)\, \, \ud x,\ \ u\in H^1(\RN).
$$
Let us set
$$
\ti{\mathcal{N}}_m:=\left\{u\in H^1(\RN)\setminus\{0\}: \ti{L}_m'(u)=0, u(0)=\max_{x\in \RN}u(x)\right\}
$$
and
$$
\ti{S}_m:=\left\{u\in\ti{\mathcal{N}}_m: \ti{L}_m(u)=\ti{E}_m\right\},\ \ \mbox{where}\ \ \ \ti{E}_m:=\inf_{u\in \ti{\mathcal{N}}_m}\ti{L}_m(u).
$$
Then, as a corollary of \cite[Lemma 2.16]{Figu}, we have the following result:
\bl\lab{lemma1}
Assume $(M1)$ if $N=2$ and $(M1)$-$(M5)$ if $N\ge3$. Then $S_m\not=\phi$ if $\ti{S}_m\not=\phi$. Moreover, there exists a injective mapping $T: \ti{S}_m\longrightarrow S_m$. In particular, $T$ is bijective for $N=2$.
\el
\br\lab{r2.1}
In \cite{Figu}, Lemma \ref{lemma1} is introduced in the subcritical case. It is easy to check that the proof does not depend on the growth of the nonlinearity $f(s)$ at infinity.
\er
Assuming that $\ti{S}_m\not=\phi$, as in \cite{Figu}, the mapping $T$ is given as follows:
\begin{itemize}
\item[(i)] If $N=2$ then $T: \ti{S}_m\longrightarrow S_m$ is given by
$$
(Tu)(\cdot):=u\left(\left(M(\|\na u\|_2^2)\right)^{-1/2}\cdot\right),\ u\in\ti{S}_m.
$$
\item[(ii)]  If $N\ge3$, $T: \ti{S}_m\longrightarrow S_m$ is defined by
$$
(Tu)(\cdot):=u(\cdot/t_u), u\in\ti{S}_m,
$$
where
$
t_u:=\inf\left\{t>0: t^2=M(t^{N-2}\|\na u\|_2^2))\right\}.
$
\end{itemize}
\bl\lab{lemma}
Assuming that $\ti{S}_m\not=\phi$ for $N\ge2$, then $S_m\not=\phi$. Moreover, for any $v\in S_m$, there exists $u\in\ti{S}_m$ such that $v(\cdot)=u(\cdot/h_v)$, where
\begin{align*}
h_v=\left\{
\begin{array}{cc}
  \sqrt{M(\|\na u\|_2^2)}\,\,&\,\mbox{if}\,\, N=2; \\
  \sqrt{M(\|\na v\|_2^2)}\,\,&\,\mbox{if}\,\, N\ge3.
\end{array}
\right.
\end{align*}
\el
\bp
By the definition of $T$, we know $S_m\not=\phi$ if $\ti{S}_m\not=\phi$. Let $v\in S_m$, then $-M(\|\na v\|_2^2)\DD v=f(v)-mv, v\in H^1(\RN)$ and $L_m(v)=E_m$. If $N=2$, then by Lemma \ref{lemma1}, $u=T^{-1}v\in\ti{S}_m$ and $v(\cdot)=u(\cdot/h)$, $h=\sqrt{M(\|\na u\|_2^2)}$.

If $N\ge3$, let
$$
u(\cdot):=v(h\cdot),\,\,\mbox{where}\,\, h=\sqrt{M(\|\na v\|_2^2)},
$$
then
$$
-\DD u(\cdot)=-h^2\DD v(h\cdot)=f(v(h\cdot))-mv(h\cdot),
$$
i. e., $-\DD u+mu=f(u)$ in $\RN$. In the following, we show that $u\in\ti{S}_m$. It suffices to show that $\ti{L}_m(u)=\ti{E}_m$. By the Pohozaev's identity,
$$
\ti{L}_m(u)=\frac{1}{N}\left[\frac{M(\|\na v\|_2^2)}{(\|\na v\|_2^2)^\frac{2}{N-2}}\right]^{\frac{2-N}{2}}.
$$
On the other hand, let $\ti{u}\in\ti{S}_m$, then $\ti{v}:=T\ti{u}=\ti{u}(\cdot/t_{\ti{u}})\in S_m$, where $t_{\ti{u}}$ is given above. By the Pohozaev's identity, we know
\be\lab{mv}\left\{
\begin{array}{c}
  L_m(\ti{v})=\frac{1}{2}[\hat{M}(\|\na\ti{v}\|_2^2)-(1-\frac{2}{N})M(\|\na\ti{v}\|_2^2)\|\na\ti{v}\|_2^2]=E_m \\
  L_m(v)=\frac{1}{2}[\hat{M}(\|\na v\|_2^2)-(1-\frac{2}{N})M(\|\na v\|_2^2)\|\na v\|_2^2]=E_m.
\end{array}
\right.
\ee
Then by the proof of \cite[Lemma 2.17]{Figu} and $(M5)$, it is easy to know that if for some $0\le t_1<t_2$, $$\hat{M}(t_1)-(1-\frac{2}{N})M(t_1)t_1=\hat{M}(t_1)-(1-\frac{2}{N})M(t_1)t_1,$$
then
$$
\frac{M(t_1)}{t_1^\frac{2}{N-2}}=\frac{M(t_2)}{t_2^\frac{2}{N-2}}.
$$
It follows from that \re{mv} that
$$
\ti{L}_m(u)=\frac{1}{N}\left[\frac{M(\|\na \ti{v}\|_2^2)}{(\|\na \ti{v}\|_2^2)^\frac{2}{N-2}}\right]^{\frac{2-N}{2}}=\frac{1}{N}\|\na \ti{u}\|_2^2=\ti{E}_m.
$$
Thus, $u\in\ti{S}_m$ and $v(\cdot)=u(\cdot/h)$. The proof is completed.
\ep

\bl\lab{lemma2}
Assume that $\ti{S}_m\not=\phi$ for $N\ge2$. Then there exits $C,c>0$(independent of $v$) such that $c\le h_v\le C$ for all $v\in S_m$, where $h_v$ is given in Lemma \ref{lemma}.
\el
\bp If $N=2$, take any $v\in S_m$, then $h_v=\sqrt{M(\|\na u\|_2^2)}$ for some $u\in\ti{S}_m$. By the Pohozaev's identity, $\ti{L}_m(u)=\frac{1}{2}\|\na u\|_2^2=\ti{E}_m$. Then $h_v=\sqrt{M(2\ti{E}_m)}$ for any $v\in S_m$. The desired result follows from $(M1)$. If $N\ge3$, take any $v\in S_m$, then $h_v=\sqrt{M(\|\na v\|_2^2)}$. By $(M1)$, $h_v\ge m_0$. On the other hand, by the Pohozaev's identity,
$$
L_m(v)=\frac{1}{2}[\hat{M}(\|\na v\|_2^2)-(1-\frac{2}{N})M(\|\na v\|_2^2)\|\na v\|_2^2]=E_m,\,\,\mbox{for all}\,\, v\in S_m.
$$
Then by $(M2)$, $\sup_{v\in S_m}h_v<\iy$. The proof is completed.
\ep

Now, we summarize some results on $\ti{S}_m$, whose proof can be found in \cite{Zhang-Zou,Byeon-Zhang-Zou,JM1}. Thanks to Lemma \ref{lemma} and Lemma \ref{lemma2}, $S_m$ has similar properties below, which will be used in the proof of Theorem \ref{Theorem 1}.

\bl \rm{(see \cite{Zhang-Zou,Byeon-Zhang-Zou,JM1})}\lab{lemma3}\ If $(F1)$-$(F3)$ hold then $\ti{S}_m\not=\phi$. Moreover,
\begin{itemize}
\item [$(i)$] any $U\in \ti{S}_m$ is such that $U\in C^2(\RN)\cap L^\iy(\RN)$ and is radially symmetric;
\item [$(ii)$] $\ti{S}_m$ is compact in $H^1(\RN)$;
\item [$(iii)$] $0<\inf\{\|U\|_{\infty}:U\in \ti{S}_m\}\le \sup\{\|U\|_{\infty}:U\in \ti{S}_m\}<\iy$;
\item [$(iv)$] there exist constants $C,c>0$ independent of $U\in \ti{S}_m$ such that
$$
|D^{\al}U(x)|\le C\exp(-c|x|), \,x\in \RN\ \ \mbox{for\  $|\al|=0,1$}.
$$
\end{itemize}
\el
\bp
For convenience of the reader, we provide some details here.
\vskip0.1in
{\bf Existence of ground state solutions: $\ti{S}_m\not=\phi.$}
\vskip0.1in
It is well known that \eqref{q51} possesses a ground state solution by means of the following constrained minimization problem
\begin{equation}\lab{cmp}
\begin{cases}
A:=\inf\left\{T_0(u): G(u)=1,u\in H^1(\RN)\setminus \{0\}\right\},\ \ \mbox{if}\ \ N\ge3,
\vspace{0.2cm}
\\
\displaystyle A:=\inf\left\{T_0(u): G(u)=0,u\in H^1(\R^2)\setminus \{0\}\right\},\ \ \mbox{if}\ \ N=2,
\end{cases}
\end{equation}
where $$T_0(u)=\frac{1}{2}\int_{\RN}|\nabla u|^2\, \ud x,\ G(u)=\int_{\RN}(F(u)-\frac{m}{2}u^2)\, \ud x.$$
If problem \re{cmp} admits a minimizer $u$, then there exists some $\si>0$ such that $u(\cdot/\si)$ is indeed a ground state solution of \eqref{q51}(see \cite{souto,Jean}).

In the following, we show that $A$ can be achieved.
\vskip0.1in
{\bf Case 1. $N\ge3$.} With $(F1)$, $(F2)$ and $(i)$ (or $(ii)$) of $(F3)$, Zhang and Zou \cite{Zhang-Zou} proved that $A$ can be achieved. If we assume $(F1)$, $(F2)$ and $(iii)$ of $(F3)$, the proof can be done by using a similar argument to that in \cite{Zhang-Zou}. Indeed, as can be seen in \cite{souto}, if one defines the mountain pass value
\be\lab{mpv}
b:=\inf_{\g\in\G}\max_{t\in[0,1]}\ti{L}_m(\g(t)),
\ee
where $\G:=\{\g\in C([0,1],H^1(\RN)): \g(0)=0, \ti{L}_m(\g(t))<0\}$, then
$$
b=\frac{1}{N}\left(\frac{N-2}{2N}\right)^{(N-2)/2}(2A)^{N/2}.
$$
By $(iii)$ of $(F3)$ we know that $b<\frac{1}{N}S^{N/2}$ for $\la>0$ large enough, where $S$ is the best Sobolev's embedding constant of $D^{1,2}(\RN)\hookrightarrow L^{2^\ast}(\RN)$. Then $0<A<\frac{1}{2}\big(2^{\ast}\big)^{\frac{N-2}{N}}S$ for $\la>0$ large enough. And, by following the argument in \cite{Zhang-Zou}, it is easy to show that $A$ is achieved.
\vskip0.1in
{\bf Case 2. $N=2$.} As can be seen in \cite{souto, Ruf}, in order to  the existence of a minimizer for $A$, it suffices to prove $A<1/2$. By \cite{souto}, we know that $A\le c$, where
$$
c:=\inf_{u\in H^1(\R^2)/\{0\}}\max_{t\ge0}\ti{L}_m(tu).
$$
In the following, we use the argument of Adimurthi \cite{Adimurthi} (see also \cite{FMR,Ruf,JFZ}) to construct a function $w\in H^1(\R^2)\setminus\{0\}$ such that $\max_{t\ge0}\ti{L}_m(tw)<1/2$, which implies that $b<1/2$. The proof is standard. Again, for convenience of the reader, we give the details. By $(F3)$, choosing some fixed $r>0$ such that
\begin{equation}
\label{ChoiceOfBeta}
\beta_0 >\frac{e^{\frac 1 2 \: r^2 m}}{\pi  r^2},
\end{equation}
we consider the Moser sequence of functions
$$
\tilde w_n(x):= \frac 1{\:\sqrt{2 \pi} \:}
\begin{cases}
\sqrt{\log n}, & \text{if } |x| \leq \frac r n;
\vspace{0.2cm}
\\
\displaystyle{\frac{\log \frac r{|x|}}{\:\sqrt{\log n}\:}}, & \text{if } \frac r n \leq |x| \leq r;
\vspace{0.2cm}
\\
0, & \text{if } |x| \geq r.
\end{cases}
$$
It is well known that $\|\na\tilde w_n\|_2=1$ and $\|\tilde w_n\|_2^2=r^2/(4\log{n})+o(r^2/\log{n})$. Let
$$
\|\tilde w_n\|^2:=\|\na\tilde w_n\|_2^2+m\|\tilde w_n\|_2^2=1 + \frac{d_n(r)}{\log n}m,
$$
where
$d_n(r):= r^2/4 + o_n(1) \quad \text{ and } \quad o_n(1) \to 0 \text{ as } n \to + \infty$. Setting
$w_n:=\tilde w_n/\|\tilde w_n\|$ then, for $n$ large enough,
\begin{equation}
\label{EstOfwn}
(w_n)^2(x)\geq \frac 1{2 \pi} \: \Bigl( \log n -  \: d_n(r)  m \Bigr) \quad \text{for } |x| \leq \frac r n.
\end{equation}
Now, we prove that there exists some $n\in\mathbb{N}$ such that $\max_{t\ge0}\ti{L}_m(tw_n)<1/2$. Assume,  on the contrary, that
$$
\max_{t \geq 0} \ti{L}_{m}(t w_n) \geq \frac{1}{2},\, \mbox{for all}\, n \in \mathbb N.
$$
As a consequence of $(F3)$, for any $\varepsilon >0$ there exists $R_\varepsilon>0$ such that
\begin{equation}
\label{consf3}
sf(s) \geq (\beta_0 - \varepsilon) e^{4\pi s^2}, \quad \forall s \geq R_\varepsilon.
\end{equation}
Then it is easy to see that $\ti{L}_{m}(t w_n)\rg -\iy$ as $t\rg\iy$. And, by our assumption, there exists $t_n >0$ such that
\be\lab{mi}
\ti{L}_{m}(t_n w_n)= \max_{t \geq 0} \ti{L}_{m}(t w_n)\ge1/2.
\ee
Noting that $\|w_n\|=1$ and $f(s)\ge0$ for all $s\ge0$, we have
$$
\frac{1}{2}t_n^2\ge\frac{1}{2}t_n^2-\int_{\R^2}F(tw_n)\ge1/2,
$$
which implies that $t_n\ge1$.
\vskip0.1in
Next, we claim that $\lim\limits_{n\rg\iy}t_n=1$. Note that
\begin{equation}
\label{tn2}
t_n^2= \int_{\R^2} f(t_n w_n) t_n w_n\, \ud x
\end{equation}
and $$t_nw_n=\frac{t_n}{\| \ti w_n\|}\frac{\sqrt{\log{n}}}{\sqrt{2\pi}}\rg+\iy\ \ \mbox{as}\ \ n\rg\iy,\ \ x\in B_{r/n},$$
for $n$ large enough. Using \re{EstOfwn}, \re{consf3} and \re{tn2}, we get for $n$ large enough that
\begin{align*}
t_n^2 &\geq (\beta_0- \varepsilon) \int_{B_{r/n}}e^{4\pi (t_n w_n)^2} \, \ud x\\
&\geq \pi r^2(\beta_0 - \varepsilon)\: e^{2 t_n^2[\log n - d_n(r)m]-2\log{n}},
\end{align*}
which implies that $\{t_n\}$ is bounded and also $\limsup\limits_{n\rg\iy}t_n\le1$. Thus, $\lim\limits_{n\rg\iy}t_n=1$.

Noting that $w_n \to 0 $ a.e. in $\R^2$, Lebesgue's dominated convergence theorem yields (as $n\rg\iy$):
$$
\int_{\{t_nw_n < R_\varepsilon\}} f(t_n w_n) t_n w_n \, \ud x \rg0 \quad \text{ and } \quad \int_{\{t_nw_n < R_\varepsilon\}} e^{4\pi (t_n w_n)^2} \, \ud x\rg\pi r^2.
$$
Then, it follows from \eqref{consf3} and \eqref{tn2} that
{\allowdisplaybreaks
\begin{equation}\lab{bh}
\begin{split}
t_n^2 & =\int_{B_r} f(t_n w_n) t_n w_n \, \ud x
\\
& \geq (\beta_0- \varepsilon) \int_{B_r} e^{4\pi (t_n w_n)^2} \, \ud x + \int_{\{t_nw_n < R_\varepsilon\}} f(t_n w_n) t_n w_n \, \ud x
\\ &\ \ \ \  - (\beta_0- \varepsilon) \int_{\{t_nw_n < R_\varepsilon\}} e^{4\pi (t_n w_n)^2} \, \ud x
\\ &\ge(\beta_0- \varepsilon)\Big[\int_{B_r} e^{4\pi (w_n)^2} \, \ud x-\pi r^2\Big],
\end{split}
\end{equation}
}%
for $n$ large enough. Also, it follows from \eqref{EstOfwn} that
$$
\liminf_{n\rg\iy}\int_{B_{r/n}}  e^{4\pi (w_n)^2} \, \ud x \geq \pi r^2e^{-m r^2/2}
$$
for $n$ large. On the other hand, using the change of variable $s=r e^{-\| \ti w_n\| \sqrt{\log n} \: t}$, we have
\begin{align*}
\int_{B_r \setminus B_{r/n}} e^{4 \pi (w_n)^2} \, \ud x & = 2 \pi r^2 \| \ti w_n\| \sqrt{\log n} \: \int_{0}^{\frac{\sqrt{\log n}}{\| \ti w_n\|}} \: e^{2(\:t^2 - \| \ti w_n\| \sqrt{\log n} \: t \:)} \, \ud t
\\
& \geq 2 \pi r^2 \| \ti w_n\| \sqrt{\log n} \: \int_{0}^{\frac{\sqrt{\log n}}{\| \ti w_n\|}} \: e^{-2\| \ti w_n\| \sqrt{\log n} \: t} \, \ud t \\
& = \pi r^2 \bigl( 1 - e^{-2 \log n}\bigr).
\end{align*}
So, by \re{bh} we have
$$
1= \lim_{n \to + \infty} t_n^2 \geq (\beta_0 - \varepsilon) \pi r^2e^{-m r^2/2}.
$$
Since $\e$ is arbitrary, we obtain
$$\beta_0 \leq \frac{e^{\frac 1 2 \: r^2m}}{\pi r^2},$$
which contradicts \eqref{ChoiceOfBeta}. Hence, $\max_{t\ge0}\ti{L}_m(tw_n)<1/2$ for some $n$, which implies that $A<1/2$. Therefore $A$ can be achieved.
\vskip0.1in
{\bf Regularity of ground state solutions.}
\vskip0.1in
In the case $N\ge3$, the properties $(i)$-$(iv)$ were given in \cite[Proposition 2.1]{Byeon-Zhang-Zou}. For $N=2$, we refer to \cite{JZ1}. Once again, for the convenience of the reader, we give a sketch of the proof in the case $N=2$.

{\bf Step 1.} For any $U\in\ti{S}_m$ we claim that $U\in L^\iy(\R^2)$.

Indeed, by the Trudinger-Moser inequality (see \cite{JM1}), $f(U)\in L_{loc}^2(\R^2)$, which implies by {\it interior $H^2$-regularity} (see \cite{Evans}) that $U\in H_{loc}^2(B_r)$. Moreover, for each open set $\om\subset\subset B_r$ with $\pa\om\in C^1$,
\be\lab{v1}
\|U\|_{H^2(\om)}\le C\left(\|f(U)\|_{L^2(B_r)}+\|U\|_{L^2(B_r)}\right),
\ee
where $C$ depends only on $\om, r$. By the Sobolev's embedding theorem, $U\in C^{0,\g}(\overline{\om})$ for some $\g\in(0,1)$ and there exists $c$ (independent of $U$) such that
\be\lab{v2}\|U\|_{C^{0,\g}(\overline{\om})}\le c\|U\|_{H^2(\om)}.
\ee
Now, we show that $\lim_{|x|\rg\iy}U(x)=0$. Suppose on the contrary that there exists $\{x_j\}\subset\R^2$ with $|x_j|\rg\iy$ as $j\rg\iy$ and $\liminf_{j\rg\iy}U(x_j)>0$. Let $v_j(x)=U(x+x_j)$, then
\be\lab{v3}
-\DD v_j+mv_j=f(v_j), v_j\in H^1(\R^2).
\ee
Assume that $v_j\rg v$ weakly in $H^1(\R^2)$. Then, by elliptic estimates we have $v\not\equiv 0$. However, for any fixed $R>0$,
\begin{align*}
\int_{\R^2}U^2&\ge\liminf_{j\rg\iy}\left(\int_{B_R(0)}U^2+\int_{B_R(x_j)}U^2\right)\\
&=\int_{B_R(0)}U^2+\liminf_{j\rg\iy}\int_{B_R(0)}v_j^2\\
&=\int_{B_R(0)}U^2+\int_{B_R(0)}v^2\\
&\rg\int_{\R^2}U^2+\int_{\R^2}v^2,\ \ \mbox{as}\ \ R\rg\iy,
\end{align*}
which is a contradiction. Thus, $U(x)\rg 0$ as $|x|\rg\iy$. Noting that $U\in C(\R^2)$, we have $U\in L^\iy(\R^2)$.

\vskip0.1in
\noindent{\bf Step 2.} For any $U\in\ti{S}_m$ we claim that $U$ is radially symmetric, which implies that $U\in C^2(\R^2)$.

Indeed, let us consider the constrained minimization problem \re{cmp} for $N=2$. For any minimizer $u$ of \eqref{cmp}, as we can see in \cite{be}, there exists $\theta>0$ such that
$$
\int_{\R^2}\na u\na \varphi=\theta\int_{\R^2}(f(u)-\frac{m}{2}u)\varphi,\ \forall \varphi\in H^1(\R^2),
$$
namely, $u$ satisfies
\be\lab{cm1}
-\DD u+\theta mu=\theta f(u),\ u\in H^1(\R^2).
\ee
Similarly to above, $u\in C(\R^2)\cap L^{\iy}(\R^2)$ and $u(x)\rg 0$ as $|x|\rg\iy$. By $C^\al$-regularity theory (see \cite[Theorem 10.1.2]{Jost}), $u\in C^{1,\al}(\R^2)$ for some $\al\in(0,1)$. Moreover, for any solution $u$ of \eqref{cm1}, $u\in C^{1,\al}(\R^2)$ and $u(x)\rg 0$ as $|x|\rg\iy$. By a classical comparison argument, $u$ decays exponentially at infinity. Then by Pohozaev's identity, $u$ satisfies $G(u)=0$. By $(F1)$, $F(s)-\frac{m}{2}s^2<0$ for small $|s|>0$. Therefore, by \cite[Proposition 4]{Maris} we know that $U$ is radially symmetric.

\vskip0.1in

\noindent{\bf Step 3.} We claim that $\ti{S}_m$ is compact in $H^1(\R^2)$.

Indeed, by first adopting some ideas in \cite{BJT}, we can prove that $\ti{S}_m$ is bounded in $H^1(\R^2)$,  so obviously, $\{\|\na U\|_{L^2(\R^2)}^2|U\in\ti{S}_m\}$ is bounded. Now we claim that $\{\|U\|_{L^2}^2|U\in\ti{S}_m\}$ is bounded. Otherwise, there exists $\{U_j\}\subset\ti{S}_m$ such that $\la_j=\|U_j\|_{L^2}\rg\iy$ as $j\rg\iy$. Letting $\wti{U}_j(x)=U_j(\la_j x)$, then $\widetilde{U}_j$ satisfies $\|\widetilde{U}_j\|_{L^2}=1$, $\|\na \widetilde{U}_j\|_{L^2}^2=2\ti{E}_m$ and
\be\lab{eq2}
-\la_j^{-2}\DD \widetilde{U}_j+m\widetilde{U}_j=f(\widetilde{U}_j)\ \mbox{in}\ \R^2.
\ee
Therefore, by $(F1)$ as in \cite{BJT}, we can assume that $\wti{U}_j\rg 0\in H_{rad}^1(\R^2)$ weakly in $H^1(\R^2)$. Noting that $\ti{E}_m=A<\frac{1}{2}$ and using a similar argument in \cite[Lemma 5.1]{souto}, it follows that $\int_{\R^2}\wti{U}_jf(\wti{U}_j)\rg 0$ as $j\rg\iy$. Thus, by \eqref{eq2}, we get $\|\wti{U}_j\|_2\rg0$ as $j\rg\iy$, which is a contradiction. Therefore, $\ti{S}_m$ is bounded in $H^1(\R^2)$. Secondly, assuming $\{u_n\}\subset \ti{S}_m$ and $u_n\rg u$ weakly in $H^1(\R^2)$, we prove that $u\in\ti{S}_m$ and, up to a subsequence, $u_n\rg u$ strongly in $H^1(\R^2)$. Obviously, it follows from \cite[Lemma 5.1]{souto} that $\int_{\R^2}F(u_n)\rg\int_{\R^2}F(u)$. Then, from $(F1)$-$(F2)$ and $0<A<1/2$, we get that $u\not\equiv0$. Noting that $u$ is a weak solution of \eqref{q51} one has $\ti{L}_m(u)\ge\ti{E}_m$. On the other hand, by Fatou's Lemma, $\ti{L}_m(u)\le \ti{E}_m$. It follows that $u\in \ti{S}_m$ and $u_n\rg u$ strongly in $H^1(\R^2)$. Therefore, $\ti{S}_m$ is compact in $H^1(\R^2)$.

\vskip0.1in

\noindent{\bf Step 4.} The property $\inf\{\|u\|_{\infty} : \ u\in \ti{S}_m\}>0$ is obvious since $\lim_{t \to 0}f(t)/t = 0.$  Noting that $\ti{S}_m$ is compact in $H^1(\R^2)$, in order to prove that $\sup\{\|u\|_{\infty} : \ u\in \ti{S}_m\}<\iy$, it suffices to prove that for any $\{u_n\}\subset \ti{S}_m$ with $u_n\rg u\in \ti{S}_m$ strongly in $H^1(\R^2)$, it holds that $\sup_n\|u_n\|_{\infty}<\iy$. First, by $(F1)$-$(F2)$, there exist $C>0$ and $\beta>4\pi$ such that $0<f(t)\le\frac{m}{2}t,\ t\in(0,1)$ and $0<f(t)\le C(\exp(\beta t^2)-1)$ for $t\ge 1$. Since $u\in L^{\iy}(\R^2)$ and $u_n\rg u$ strongly in $H^1(\R^2)$, we have by the Trudinger-Moser inequality (see \cite{JM}) that
\be\lab{yf1}
\lim_{n\rg\iy}\int_{\R^2}|\exp(2\beta u_n^2)-\exp(2\beta u^2)|=0.
\ee
Secondly, we claim that
\be\lab{ii}
\sup_n\|f(u_n)\|_2<\iy.
\ee
Letting
$$
A_n:=\{x\in\R^2|u_n(x)\le 1\},\quad B_n:=\{x\in\R^2|u_n(x)>1\},
$$
then
\begin{align*}
\int_{\R^2}|f(u_n)|^2&=\int_{A_n}|f(u_n)|^2+\int_{B_n}|f(u_n)|^2\\
&\le\int_{\R^2}\frac{m^2}{4}|u_n|^2+C\int_{\R^2}[\exp(2\beta u_n^2)-1].
\end{align*}
It follows by the Trudinger-Moser inequality (\cite{JM}) and by \eqref{yf1} that \eqref{ii} is true. Similarly, by {\it interior $H^2$-regularity} (see \cite{Evans}), we have that
\be\lab{v11}
\|u_n\|_{H^2(B_1)}\le C\left(\|f(u_n)\|_{L^2(B_2)}+\|u_n\|_{L^2(B_2)}\right),
\ee
where $C$ is independent of $n$. On the other hand, by the Sobolev's embedding theorem,
\be\lab{v21}
\|u_n\|_{C^{0,\g}(\overline{B_1})}\le c\|u_n\|_{H^2(B_1)}
\ee
for some $\g\in(0,1)$, where $c$ is independent of $n$. Therefore, it follows from \eqref{ii}-\eqref{v21} that $\sup_n\|u_n\|_{C^{0,\g}(\overline{B_1})}<\iy,$ which implies that, up to a subsequence, $u_n\rg u$ uniformly in $\overline{B_1}$. Thus, since $u\in L^{\iy}(\R^2)$, we get that $\sup_n\|u_n\|_{L^{\iy}(\overline{B_1})}<\iy$. Therefore, $\sup_n\|u_n\|_{L^{\iy}(\R^2)}<\iy$.

\vskip0.1in

\noindent{\bf Step 5.} By the radial lemma \cite{Strauss}, $u_n(x)\rg 0$ as $|x|\rg\iy$ uniformly w.r.t. $n$. By a classical comparison principle, $\sup_n\|u_n\|_{L^{\iy}(\R^2)}<\iy$ and there exist $c,C>0$ such that
$$
U(x)+|\na U(x)|\le C\exp(-c|x|), \,x\in \R^2,
$$
for any $U\in\ti{S}_m$. The proof is complete.
\ep

\vskip0.1in

\s{Proof of Theorem \ref{Theorem 1}}

\renewcommand{\theequation}{3.\arabic{equation}}
By Proposition \ref{prop1}, there exists $\kappa>0$ such that
\be\lab{ee}\sup_{U\in S_m}\|U\|_\iy=\sup_{U\in \ti{S}_m}\|U\|_\iy<\kappa.\ee
For any fixed $k>\max_{t\in[0,\kappa]}f(t)$, define $f_k(t)=\min\{f(t),k\},\, t\in\R.$
Now, we consider the truncated problem

\be\lab{q6}
-\e^2M\left(\e^{2-N}\int_{\RN}|\na v|^2\, \ud x\right)\Delta v+V(x)v=f_k(v)\ \mbox{in}\ \RN.
\ee
In the following we prove that, for small $\e>0$, there exists a positive solution $v_\e$ of \re{q6} satisfying the properties $(i)$-$(ii)$ in Theorem \ref{Theorem 1}. Obviously, $v_\e$ is a solution of the original problem \re{q2} if $\|v_\e\|_{\iy}<\kappa$.

We consider the limiting problem of \re{q6}
\be\lab{q7}
-M(\|\na u\|_2^2)\DD u+mu=f_k(u),\ \ u\in H^1(\RN),
\ee
whose energy functional is given by
$$
\ti{L}_m^k(u)=\frac{1}{2}\hat{M}(\|\na u\|_2^2)+\frac{m}{2}\int_{\RN}u^2-\int_{\RN}F_k(u),\ \ u\in H^1(\RN),
$$
where $F_k(s)=\int_0^sf_k(t)\, \ud t$.

\bl\lab{l5}
With the same assumptions in Theorem \ref{Theorem 1}, the limit problem \eqref{q7} admits one positive ground state solution, which is radially symmetric.
\el
\bp
By the definition of $f_k$, it is easy to check that $f_k$ satisfies $(f_1)$-$(f_2)$ in Theorem {\bf A}. Let $U\in \ti{S}_m$. By Pohozaev's identity,
\begin{equation*}
\begin{cases}
\displaystyle\int_{\RN}(F(U)-\frac{m}{2}U^2)\, \ud x=\frac{N-2}{2N}\int_{\RN}|\na U|^2\, \ud x\ \ \mbox{if}\ \ N\ge3,
\vspace{0.2cm}
\\
\displaystyle\int_{\R^2}(F(U)-\frac{m}{2}U^2)\, \ud x=0\, \ud x\ \ \mbox{if}\ \ N=2,
\end{cases}
\end{equation*}
we get that $\int_{\RN}(F(U)-\frac{m}{2}U^2)\, \ud x\ge0$ for $N\ge2$. If $F(U(x))\le\frac{m}{2}U^2(x)$ for all $x\in\RN$, then $F(U(x))/U^2(x)\equiv m/2>0$ for all $x\in\RN$. Recalling that $U(x)\rg0$ as $|x|\rg\iy$, by $(F1)$ we get that $F(U(x))/U^2(x)\rg0$ as $|x|\rg\iy$, which is a contradiction. So there exists $x_0\in\RN$ such that $F(U(x_0))>\frac{m}{2}U^2(x_0)$. Noting that $|U(x_0)|<\kappa$, we have $F_k(U(x))\equiv F(U(x))$ for all $x\in\RN$. Then, letting $\xi=U(x_0)>0$, we have $F_k(\xi)>\frac{m}{2}\xi^2$. Namely, $f_k$ satisfies $(f_3)$ in Theorem {\bf A}. Therefore, it follows from \cite{Lions,be} that the problem
\be\lab{q77}
-\DD u+mu=f_k(u),\ \ u\in H^1(\RN)
\ee
admits a radially symmetric ground state solution. By \cite[Lemma 2.16]{Figu}, the proof is finished.
\ep

Let $S_m^k$ be the set of positive ground state solutions $U$ of \eqref{q7} satisfying $U(0)=\max_{x\in \R^2}U(x)$. Then by Lemma \ref{l5} $S_m^k\not=\phi$.
\bl\lab{l6}
For $k>\max_{t\in[0,\kappa]}f(t)$, we have
$$
S_m^k=S_m.
$$
\el
\bp
By Lemma \ref{lemma1} and \ref{lemma}, it suffices to prove $\ti{S}_m^k=\ti{S}_m$. Denote by $\ti{E}_m^k$ the least energy of \re{q77} and by $\ti{S}_m^k$ the set of positive ground state solutions $U$ of \eqref{q77} with $U(0)=\max_{x\in \RN}U(x)$. It follows from \cite{Jean} and \cite{Zhang-Zou} that $\ti{E}_m^k$ and $E_m^k$ coincide with the mountain pass values respectively. Noting that $f_k(t)\le f(t)$ for any $t$, we have $\ti{E}_m^k\ge \ti{E}_m$. By the definition of $f_k$, $f_k(u)=f(u)$ for any $u\in\ti{S}_m$. Then $u$ is a nontrivial solution of \re{q77}, which implies $\ti{E}_m^k\le\ti{E}_m$. Therefore
$$
\ti{E}_m^k=\ti{E}_m \ \mbox{for}\ k>\max_{t\in[0,\kappa]}f(t).
$$
Obviously, $\ti{S}_m\subset\ti{S}_m^k$ for $k>\max_{t\in[0,\kappa]}f(t)$. In what follows, we show that $\ti{S}_m^k\subset\ti{S}_m$ for $k>\max_{t\in[0,\kappa]}f(t)$. Let
$$
G_k(u)=\int_{\RN}(F_k(u)-\frac{m}{2}|u|^2)\, \ud x.
$$
Then, by \cite{Jean}(see also \cite{souto}), we have
$$
\ti{E}_m^k=
\begin{cases}
\frac{1}{N}\left(\frac{N-2}{2N}\right)^{(N-2)/2}(2A_k)^{N/2}, & \text{if } N\ge3,
\vspace{0.2cm}
\\
\displaystyle A_k, & \text{if } N=2,
\end{cases}
$$
where
\begin{equation}\lab{m1}
\begin{cases}
A_k=\inf\left\{T_0(u): G_k(u)=1,u\in H^1(\RN)\setminus\{0\}\right\}\,\mbox{if}\, N\ge3,
\vspace{0.2cm}
\\
\displaystyle A_k=\inf\left\{T_0(u): G_k(u)=0,u\in H^1(\R^2)\setminus\{0\}\right\}\,\mbox{if}\, N=2.
\end{cases}
\end{equation}
Now, we consider the cases: $N\ge3$ and $N=2$ separately.
\vskip0.1in
{\bf Case 1: $N\ge3$}.
\vskip0.1in
For any $u\in \ti{S}_m^k$, by Pohozaev's identity, we have $\|\na u\|_2^2=N\ti{E}_m^k$ and $G_k(u)=\frac{N-2}{2}\ti{E}_m^k$. Let $v(\cdot)=u(\si\cdot)$, where $\si=\left(\frac{N-2}{N}A_k\right)^{1/2}$. Then $\|\na v\|_2^2=2A_k$ and $G_k(v)=1$. Recalling that $\ti{E}_m^k=\ti{E}_m$, we have $A_k=A$. Therefore $v$ satisfies $T_0(v)=A$ and $G(v)\ge 1$. On the other hand, note that
\be\lab{m2}
A=\inf\left\{T_0(u):G(u)=1,u\in H^1(\RN)\setminus\{0\}\right\}.
\ee
If $G(v)>1$, there exists $\theta\in(0,1)$ such that $G(v_\theta)=1$, where $v_\theta(\cdot)=v(\cdot/\theta)$. However, $T_0(v_\theta)=\theta^{N-2}A<A$, which contradicts \re{m2}. So, $G(v)=1$, which implies that $v$ is a minimizer for $A$. Therefore, as can be seen in \cite{Zhang-Zou,Jean}, there exists $\si_0>0$ such that $v_{\si_0}(\cdot)=v(\cdot/\si_0)$ is a ground state solution of \re{q51}, i.e., $v_{\si_0}\in \ti{S}_m$. By \re{ee}, $\|u\|_\iy=\|v_{\si_0}\|<\kappa$, which implies that $u$ is a ground state solution of \re{q51}, i.e., $u\in \ti{S}_m$. Thus, $\ti{S}_m^k\subset\ti{S}_m$.

\vskip0.1in
{\bf Case 2: $N=2$}.
\vskip0.1in
For any $u\in \ti{S}_m^k$, by Pohozaev's identity, we have $\|\na u\|_2^2=2\ti{E}_m^k=2A_k$ and $G_k(u)=0$. Since $\ti{E}_m^k=\ti{E}_m$, $u$ satisfies $T_0(u)=A$ and $G(u)\ge 0$. Recall that
\be\lab{m22}
A=\inf\left\{T_0(u):G(u)=0,u\in H^1(\R^2)\setminus\{0\}\right\}.
\ee
If $G(u)>0$, similarly as in \cite{souto}, there exists $\theta\in(0,1)$ such that $G(\theta u)=0$. However, $T_0(\theta u)=\theta^2A<A$, which contradicts \re{m22}. So, $G(u)=0$, which implies that $u$ is a minimizer of \re{m22}. Then, as can be seen in \cite{be,souto}, there exists $\theta_0>0$ such that $u_{\theta_0}(\cdot)=u(\cdot/\sqrt{\theta_0})$ is a ground state solution of \re{q51}, i.e., $u_{\theta_0}\in \ti{S}_m$. Thus, by \re{ee} $\|u\|_\iy=\|u_{\theta_0}\|<\kappa$, which implies $u\in\ti{S}_m$. Thus $\ti{S}_m^k\subset\ti{S}_m$.
\ep

\vskip0.1in

\noindent{\bf Completion of the proof for Theorem \ref{Theorem 1}}
\vskip0.1in
\noindent{\it Proof.} First, we consider the truncation problem \re{q6}. By the proof of Lemma \ref{l5}, $f_k$ satisfies $(f_1)$-$(f_3)$ in Theorem {\bf A}. It follows from Theorem {\bf A} that for fixed $k>\max_{t\in[0,\kappa]}f(t)$, there exists $\e_0>0$ such that \eqref{q6} admits a positive solution $v_{\e}$ for $\e\in(0,\e_0)$. Moreover, there exist $U\in S_m^k$ and a maximum point $x_\e\in\RN$ of $v_\e$, such that $\lim_{\e\rg 0}\mbox{dist}(x_\e,\mathcal{M})=0$ and $v_\e(\e\cdot+x_\e)\rg U(\cdot+z_0)$ as $\e\rg 0$ in $H^1(\RN)$, for some $z_0\in\RN$. Letting $w_\e(\cdot)=v_\e(\e\cdot+x_\e)$, then $w_\e$ satisfies
$$
-M(\|\na w_\e\|_2^2)\DD w_\e+V_\e(x+\frac{x_\e}{\e})w_\e=f_k(w_\e),\ \ w_\e\in H^1(\RN).
$$
Clearly,
$$
m_0\le\inf_{\e<\e_0}M(\|\na w_\e\|_2^2)\le\sup_{\e<\e_0}M(\|\na w_\e\|_2^2)<\iy.
$$
Since $f_k(w_\e(x))\in[0,k],\ x\in\RN$, it follows from elliptic estimates that $w_\e(\cdot)\rg U(\cdot+z_0)$ locally uniformly in $\RN$. Therefore $\|v_\e\|_\iy=w_\e(0)\rg U(z_0)$ as $\e\rg0$. By Lemma \ref{l6} we have $S_m^k=S_m$, hence $U\in S_m$. By \eqref{ee}, there exists $\e^\ast<\e_0$ such that  $\|v_\e\|_\iy<\kappa$ for $\e<\e^\ast$, which implies $f_k(v_\e(x))\equiv f(v_\e(x)),\ x\in\RN$ for $\e<\e^\ast$. Therefore, $v_\e$ is a positive solution of the original problem \eqref{q2}. The proof of Theorem 1.1 is complete.
\qed

\vskip0.2in

\noindent{\bf Acknowledgement.}  The authors would like to express their sincere gratitude to the anonymous referee for his/her valuable suggestions and comments.

\end{document}